\documentclass[preprint,a4paper,12pt]{elsarticle}

\usepackage{amsmath}
\usepackage{amsfonts}
\usepackage{amssymb}
\usepackage{amsthm}

\usepackage{float}
\usepackage{graphicx}
\usepackage{multirow}
\usepackage{subcaption}

\theoremstyle{definition}

\usepackage{xcolor}
\usepackage[colorlinks=true, allcolors=blue]{hyperref}

\newcommand{\bracketround}[1]{\left(#1\right)}
\newcommand{\bracketcurly}[1]{\left\{#1\right\}}

\newcommand{\abs}[1]{\left|#1\right|}
\newcommand{\norm}[1]{\left\lVert#1\right\rVert}
\newcommand{\innerproduct}[1]{\left\langle#1\right\rangle}

\DeclareMathOperator{\sech}{sech}

\journal{Applied Mathematics and Computation}

\begin{document}

\begin{frontmatter}

\title{Neural networks for bifurcation and linear stability analysis of steady states in partial differential equations}

\author[inst1]{Muhammad Luthfi Shahab}

\author[inst1]{Hadi Susanto\corref{mycorrespondingauthor}}
\cortext[mycorrespondingauthor]{Corresponding author}
\ead{hadi.susanto@ku.ac.ae}

\affiliation[inst1]{
    organization={Department of Mathematics, Khalifa University of Science \& Technology},
    city={Abu Dhabi},
    postcode={PO Box 127788}, 
    country={United Arab Emirates}
}

\begin{abstract}
This research introduces an extended application of neural networks for solving nonlinear partial differential equations (PDEs).
A neural network, combined with a pseudo-arclength continuation, is proposed to construct bifurcation diagrams from parameterized nonlinear PDEs. Additionally, a neural network approach is also presented for solving eigenvalue problems to analyze solution linear stability, focusing on identifying the largest eigenvalue.
The effectiveness of the proposed neural network is examined through experiments on the Bratu equation and the Burgers equation. Results from a finite difference method are also presented as comparison. Varying numbers of grid points are employed in each case to assess the behavior and accuracy of both the neural network and the finite difference method.
The experimental results demonstrate that the proposed neural network produces better solutions, generates more accurate bifurcation diagrams, has reasonable computational times, and proves effective for linear stability analysis.
\end{abstract}

\begin{keyword}
Neural networks \sep Continuation \sep Bifurcation \sep Linear stability \sep Nonlinear partial differential equations \sep Bratu equation \sep Burgers equation
\MSC 65N25 \sep 65N75 \sep 65P30 \sep 65P40
\end{keyword}

\end{frontmatter}

\section{Introduction} \label{intro}

Machine learning or neural networks for approximating solutions of partial differential equations (PDEs) began in the 1990s.
Some of the earliest studies on this approach were due to Dissanayake et al.\ \cite{dissanayake1994neural} and Lagaris et al.\ \cite{lagaris1998artificial}. They used neural networks, collocation points, and finite difference formulas to transform PDEs into unconstrained optimization problems and solved them using a quasi-Newton method.
The well-known theory about universal function approximators explained by Hornik et al.\ \cite{hornik1989multilayer} supports neural networks' success.

The neural network approach offers several advantages over traditional methods like finite difference methods. 
First, neural networks provide explicit approximate solutions, allowing the value at any point to be calculated straightforwardly without interpolation \cite{lagaris1998artificial}. Additionally, neural networks support mesh-free or non-uniform grids or collocation points.
Moreover, since neural networks yield explicit approximate solutions, exact derivatives can be obtained through automatic differentiation \cite{baydin2018automatic}. Alternatively, finite difference formulas can still be used with higher accuracy for derivative calculations.
Deep neural networks are particularly adept at handling high-dimensional PDEs \cite{han2018solving, putri2024deep}. On the other hand, traditional numerical methods often struggle with high-dimensional PDEs due to the exponential increase in computational complexity and memory requirements, a challenge known as the curse of dimensionality.

Research using machine learning and neural networks to solve PDEs has significantly increased recently.
Physics-informed neural networks (PINNs) proposed by Raissi et al.\ \cite{raissi2019physics} are currently one of the most popular research in this field. It is a framework based on deep neural networks that can address direct and inverse problems related to nonlinear PDEs. They implemented the framework to deal with the Schr\"{o}dinger and the Allen-Cahn equations. However, although they explained PINNs with several new concepts, some parts also have similarities with the method proposed by Dissanayake et al.\ \cite{dissanayake1994neural}. 
PINNs have been used in various areas, including fluid mechanics \cite{cai2021physics}, power systems \cite{misyris2020physics}, optics and metamaterials \cite{jiang2022physics, chen2020physics}, finance \cite{bai2022application}, subsurface transport \cite{he2020physics}, continuum micromechanics \cite{henkes2022physics}, and supersonic flows \cite{jagtap2022physics}.
In addition, Lu et al.\ \cite{lu2021deepxde} introduced DeepXDE, a Python library designed for PINNs, offering customization and user-friendliness. Additionally, they proposed a novel residual-based adaptive refinement technique to enhance the training efficiency of PINNs.
A work proposed by Meer et al.\ \cite{van2022optimally} can be used to determine better parameters for weighted loss functions to enhance the accuracy of neural networks for solving PDEs. 
Rodriguez-Torrado et al.\ \cite{rodriguez2022physics} introduced a novel network architecture called physics-informed attention-based neural networks (PIANNs). PIANNs combine the capabilities of recurrent neural networks and attention mechanisms to capture the extreme dynamics of nonlinear hyperbolic PDEs effectively.
Yang et al.\ \cite{yang2021b} introduced a novel approach called Bayesian physics-informed neural network (B-PINN) for solving PDEs. In B-PINN, a combination of a Bayesian neural network and a PINN is employed as the prior. At the same time, the posterior can be estimated using either Hamiltonian Monte Carlo or variational inference methods.
Han et al.\ \cite{han2017deep, han2018solving} used a unique and parallel deep learning architecture for solving high-dimensional PDEs by reformulating the equations as backward stochastic differential equations \cite{pardoux2005backward}. It was then used to tackle the nonlinear Black-Scholes equation and the Hamilton–Jacobi–Bellman equation \cite{putri2024deep}. Furthermore, a review to examine the current trends in integrating physics with machine learning, outlining their existing capabilities and limitations, is given by Karniadakis et al.\ \cite{karniadakis2021physics}. The advancement of neural networks and deep learning is further facilitated by the availability of specialized hardware and software tools, like Matlab with its deep learning toolbox and TensorFlow \cite{abadi2016tensorflow} with Keras \cite{chollet2015keras,chollet2017xception}.

Despite the significant interest in applying machine learning to PDEs, researchers still underexplore some areas. One such area is the construction of bifurcation diagrams, which are crucial for understanding the behavior of solutions in parameterized PDEs. Fabiani et al.\ \cite{fabiani2021numerical} addressed this by using extreme learning machines to construct bifurcation diagrams for the Bratu equation. However, their method did not produce a complete bifurcation diagram for the two-dimensional Bratu equation. Additionally, there is research on using neural networks to effectively track the evolution of bifurcating phenomena in computational fluid dynamics problems using a Reduced Order Models method \cite{pichi2023artificial}, which differs from our proposed approach.

Another area that can be further explored is machine learning for solving eigenvalue problems and, furthermore, for linear stability analysis. While several studies have addressed eigenvalue problems, none applied their concepts to linear stability analysis. One existing research, proposed by Ben-Shaul et al.\ \cite{ben2023deep}, utilized deep learning for solving eigenvalue problems with only positive eigenvalues, which may not be applicable to linear stability analysis. Yang et al.\ \cite{yang2023neural} employed neural networks based on power methods for solving linear eigenvalue problems, which differs from our proposed approach. Other research \cite{han2020solving, zhang2022solving} discussed eigenvalue problems related to stochastic differential equations, allowing for solutions using deep learning proposed by Han et al.\ \cite{han2018solving}. However, this method is also not applicable to general eigenvalue problems.

To address this gap, we plan to combine neural networks with a pseudo-arclength continuation to construct bifurcation diagrams. Additionally, we aim to explore the application of machine learning in linear stability analysis. Our approach involves utilizing neural networks to solve the associated eigenvalue problems of the PDEs, focusing on identifying the largest eigenvalue. 
To validate our proposed method, we will conduct experiments using the one- and two-dimensional Bratu equations \cite{boyd2011one, iqbal2020efficient} and the one-dimensional viscous Burgers equation \cite{allen2013numerical}.

This paper is organized as follows.
Section\ \ref{sec2} presents the Bratu equation and the viscous Burgers equation and their eigenvalue problems.
Section\ \ref{sec3} shows the finite difference method for solving the Bratu equation and the viscous Burgers equation.
Section\ \ref{sec4} discusses the formulation of neural networks for bifurcation and linear stability analysis.
Experimental results and discussions are provided in Section\ \ref{sec5} and the conclusion is given in Section\ \ref{sec6}.

\section{Mathematical Equations} \label{sec2}

\subsection{The Bratu Equation}

The Bratu equation \cite{bratu1914equations} (also called the Liouville-Bratu-Gelfand equation or the Bratu-Gelfand equation) is a nonlinear parabolic PDE with a positive parameter $C$ and the Dirichlet boundary condition defined by
\begin{equation} \label{bratu}
	\begin{split}
		\Delta u + Ce^u & = u_t, \quad\,\,\,\, \textbf{x} \in \Omega, \\ 
		u & = 0, \qquad \textbf{x} \in \partial\Omega,
	\end{split}
\end{equation}
where $\textbf{x}=(x_1,\dots,x_d)$, $\Omega = [0,1]^d$, $t \ge 0$, $d$ is the dimension, and 
\begin{equation}
\Delta u = \sum_{i=1}^{d}\frac{\partial^2 u}{\partial x_i^2}.
\end{equation}
The analytic steady-state solution for the one-dimensional Bratu equation is given by \cite{mohsen2014simple}
\begin{equation}\label{analyticalbratu}
    u(x) = 2 \ln\bracketround{ \frac{\cosh \omega}{\cosh (\omega(1-2x))} },
\end{equation}
where $\omega$ satisfies
\begin{equation} \label{coshtheta}
    \cosh \omega = \frac{4 \omega}{\sqrt{2C}}.
\end{equation}
The analytic steady-state solution for the two-dimensional case has not been found so far.

The Bratu equation is solvable for $0< C \le C^*$ where $C^*$ is a critical point. The critical points for one- and two-dimensional Bratu equations are 3.513830719 and 6.808124423, respectively \cite{mohsen2014simple,iqbal2020numerical}. Furthermore, for one- and two-dimensional cases, the Bratu equation has two steady-state solutions (lower and upper solutions) for $0< C<C^*$, one solution for $C=C^*$, and no solution for $C>C^*$. 
Especially for higher-dimensional cases of the Bratu equation with spherical domains, 
the existence and multiplicity of the steady-state solutions have been studied extensively in \cite{bebernes2013mathematical} by transforming the Bratu equation to the equivalent second-order ordinary differential equation using the radially symmetric property.

\subsection{Eigenvalue Problem for the Linear Stability of the Bratu Equation}

In order to determine the linear stability of a steady-state solution $u(\textbf{x})$ of the Bratu equation, we consider the perturbed solution
\begin{equation} \label{u_tilde}
    \Tilde{u}(\textbf{x},t) = u(\textbf{x}) + \varepsilon w(\textbf{x},t), \qquad \abs{\varepsilon} \ll 1.
\end{equation}
Substituting Eq.\ \eqref{u_tilde} into Eq.\ \eqref{bratu} and using the Maclaurin series, we obtain the linear equation for $w(\textbf{x},t)$,
\begin{equation}
    \Delta w + C e^u w = w_t.
\end{equation}
Writing $w=e^{\lambda t} v(\textbf{x})$, then the above equation becomes the eigenvalue problem
\begin{equation} \label{eigenvalueproblem}
    \Delta v + C e^u v = \lambda v,
\end{equation}
together with the boundary condition $v=0$ on $\partial\Omega$. Here, $\lambda$ represents an eigenvalue and $v$ is its associated eigenfunction. 
A steady-state solution $u$ of the Bratu equation is linearly stable if $\max(\lambda)<0$ and linearly unstable otherwise. In the appendix, we solve the eigenvalue problem in the one-dimensional case using asymptotic analysis. 

To solve the eigenvalue problem Eq.\ \eqref{eigenvalueproblem} using an optimization approach, we write $\lambda$ using the Rayleigh quotient which, in this case, is stated as
\begin{equation} \label{Rayleighquotient}
    \lambda = \frac{\innerproduct{\Delta v + C e^u v, v}}{\innerproduct{v, v}}.
\end{equation}
The inner product is defined as
\begin{equation}
    \innerproduct{v_1, v_2} = \int_\Omega v_1 v_2 \, d\textbf{x},
\end{equation}
and it can be approximated numerically by the trapezoidal rule \cite{burden2015numerical}.

\subsection{The Viscous Burgers Equation} \label{equation_burgers1}
The one-dimensional viscous Burgers equation is a nonlinear parabolic PDE with a viscosity  parameter $\nu$ defined by \cite{allen2013numerical}
\begin{equation} \label{burgers}
	\nu \frac{\partial^2 u}{\partial x^2} - u \frac{\partial u}{\partial x} = u_t, \qquad x \in \Omega,
\end{equation}
where $\Omega = [0,1]$ and $t \ge 0$. For the purpose of numerical simulations in the next parts, we consider two types of boundary conditions. The first one is the Dirichlet boundary condition
\begin{equation} \label{burgers_BC1}
	u(0) = \rho , \qquad u(1) = 0,
\end{equation}
where $\rho>0$. The Burgers equation with this boundary condition has a unique analytic steady-state solution
\begin{equation}
	u(x) = \frac{2}{1+\exp \bracketround{ \frac{x-1}{\nu} }} - 1,
\end{equation} 
where
\begin{equation} \label{burgers_rho}
	\rho = \frac{2}{1+\exp \bracketround{ \frac{-1}{\nu} }} - 1.
\end{equation}
This case will be used to test the performance of the proposed method for different values of $\nu$.

The second boundary condition is the mixed, Neumann-Dirichlet, boundary condition given by
\begin{equation} \label{burgers_BC2}
	\frac{\partial u}{\partial x}(0) = -\varphi , \qquad u(1) = 0.
\end{equation}
The analytic steady-state solution is given by \cite{allen2013numerical}
\begin{equation}
	u(x) = \sqrt{2c} \tanh \bracketround{ \frac{\sqrt{2c}}{2\nu} (1-x) } ,
\end{equation}
where $c$ satisfies
\begin{equation} \label{theta_and_c}
	\varphi = \frac{c}{\nu} \sech^2 \bracketround{ \frac{\sqrt{2c}}{2\nu} } .
\end{equation}
This case will be used as a test case for constructing bifurcation diagrams.

The parameter $\nu$ in the Burgers equation \eqref{burgers} can be scaled to unity. However, for numerical stability illustrations, we leave it as is. The parameter combined with the boundary conditions \eqref{burgers_BC1} and \eqref{burgers_BC2} creates a boundary layer problem \cite{bender2013advanced} that, in general, requires careful numerical methods \cite{farrell2000robust}.

\subsection{Eigenvalue Problem for the Linear Stability of the Burgers Equation} \label{equation_burgers2}
To determine the linear stability of a steady-state solution $u(x)$ of the viscous Burgers equation with the mixed boundary condition, we consider the perturbed solution
\begin{equation} \label{u_tilde_burgers}
	\tilde{u}(x,t) = u(x) + \varepsilon w(x,t), \qquad \abs{\varepsilon} \ll 1.
\end{equation}
Substituting Eq.\ \eqref{u_tilde_burgers} into Eq.\ \eqref{burgers}, we obtain the linear equation for $w(x,t)$,
\begin{equation}
	\nu w_{xx} - u w_x - u_x w = w_t .
\end{equation}
Writing $w=e^{\lambda t} v(x)$, then the above equation becomes the eigenvalue problem
\begin{equation} \label{eigenvalueproblem_burgers}
	\nu v_{xx} - u v_x - u_x v = \lambda v ,
\end{equation}
together with the boundary conditions $v_x(0)=0$ and $v(1)=0$. Here, $\lambda$ represents an eigenvalue and $v$ is its associated eigenfunction.
The Rayleigh quotient for this case is stated by
\begin{equation} \label{rayleigh_burgers}
	\lambda = \frac{\innerproduct{\nu v_{xx} - u v_x - u_x v, v}}{\innerproduct{v, v}}.
\end{equation}
We also solve the eigenvalue problem \eqref{eigenvalueproblem_burgers} using asymptotic analysis in the appendix.

\section{Finite Difference Method} \label{sec3}

Although the focus of this study is the use of neural networks in solving differential equations, we will also employ the finite difference method and compare it with the results obtained by the neural network to determine the performance of the proposed method.

In the finite difference method, the continuous domain of the problem is discretized into grid points. The derivatives at the grid points in the equation are approximated using finite difference formulas, which replace the derivatives with differences between function values at nearby grid points. Replacing the derivatives with these finite difference formulas transforms the equation into a system of equations. If the system of equations is linear, it can be solved numerically using various methods, such as Gaussian elimination, Jacobi method, or Gauss-Seidel method \cite{burden2015numerical}. Newton's method is the common optimization method used for nonlinear systems.

\subsection{One-Dimensional Bratu Equation}

The domain $[0,1]$ is divided into $n$ equal subintervals of length $h=1/n$ to get steady-state solutions for the one-dimensional Bratu equation. The second derivative with respect to $x$ at a grid point $x_i=ih$ can be approximated by the 
difference formula \cite{burden2015numerical} given by
\begin{equation} \label{second_derivative}
    \frac{\partial^2 u}{\partial x^2} (x_i) \approx 
    \frac{u(x_{i+1}) - 2u(x_{i}) + u(x_{i-1})}{h^2},
\end{equation}
for $i=1,\dots,n-1$. With this formula and representing $u(x_i)$ as $u_i$, the one-dimensional Bratu equation can be transformed into a system of nonlinear equations
\begin{equation} \label{nlsfd1}
    \frac{ u_{i+1}-2u_i+u_{i-1} }{h^2} + C e^{u_i} = 0,
\end{equation}
for $i=1,\dots,n-1$ where $u_0=u_n=0$ is used to satisfy the boundary condition.
The values $u_i$, $i = 1,\dots,n-1$, are the variables that need to be sought for by Newton's method. 
Furthermore, because the number of different equations in the nonlinear system given in Eq.\ \eqref{nlsfd1} is also $n-1$, then the size of the Jacobian matrix for this problem is $(n-1) \times (n-1)$.

\subsection{Two-Dimensional Bratu Equation}

The finite difference method for the two-dimensional Bratu equation can be obtained similarly to that of the one-dimensional case. Since the domain is $[0,1]^2$, we divide both intervals for the $x$-axis and $y$-axis into $n$ equal subintervals. If the length of the subinterval is $h=1/n$, then the second derivatives with respect to $x$ and $y$ at a grid point $(x_i,y_j)=(ih,jh)$,
can be approximated by the 
difference formula \cite{burden2015numerical} given by
\begin{equation}
    \frac{\partial^2 u}{\partial x^2}(x_i,y_j) \approx 
    \frac{u(x_{i+1},y_{j}) - 2u(x_{i},y_{j}) + u(x_{i-1},y_{j})}{h^2},
\end{equation}
and
\begin{equation}
    \frac{\partial^2 u}{\partial y^2}(x_i,y_j) \approx
    \frac{u(x_{i},y_{j+1}) - 2u(x_{i},y_{j}) + u(x_{i},y_{j-1})}{h^2},
\end{equation}
for $i = 1,\dots,n-1$ and $j = 1,\dots,n-1$.
With the above formulas and representing $u(x_i,y_j)$ as $u_{i,j}$, the two-dimensional Bratu equation can be transformed into a system of nonlinear equations
\begin{equation} \label{nlsfd2}
    \frac{ u_{i+1,j} - 2u_{i,j} + u_{i-1,j} }{h^2} 
    + \frac{ u_{i,j+1} - 2u_{i,j} + u_{i,j-1} }{h^2} 
    + Ce^{u_{i,j}} = 0,
\end{equation}
for $i = 1,\dots,n-1$ and $j = 1,\dots,n-1$. The values at each grid point on the boundary are set to be zero, $u_{i,j} = 0$ for $i = 0, n$ or $j = 0, n$, to satisfy the boundary condition. The values $u_{i,j}$, $i = 1,\dots,n-1$ and $j = 1,\dots,n-1$, are the variables that need to be sought for by Newton's method. In this case, we need to optimize $(n-1)^2$ variables simultaneously. Furthermore, because the number of different equations in the nonlinear system given in Eq.\ \eqref{nlsfd2} is also $(n-1)^2$, then the size of the Jacobian matrix for this problem is $(n-1)^2\times (n-1)^2$.

\subsection{One-Dimensional Burgers Equation} \label{finite_difference_burgers}

To get steady-state solutions for the one-dimensional Burgers equation, we use processes similar to those in the one-dimensional Bratu equation. Using the 
difference formula for the first derivative \cite{burden2015numerical}
\begin{equation}
	\frac{\partial u}{\partial x} (x_i) \approx 
	\frac{u(x_{i+1}) - u(x_{i-1})}{2h},
\end{equation}
and for the second derivative in Eq.\ \eqref{second_derivative}, the one-dimensional Burgers equation can be transformed into a system of nonlinear equations
\begin{equation} \label{system_burgers}
	\nu \frac{ u_{i+1}-2u_i+u_{i-1} }{h^2} - u_i \frac{u_{i+1} - u_{i-1}}{2h} = 0,
\end{equation}
for $i=1,\dots,n-1$.
To satisfy the first boundary condition in Eq.\ \eqref{burgers_BC1}, we set $u_0=\rho$ and $u_n=0$. For the second boundary condition in Eq.\ \eqref{burgers_BC2}, we use $u_0=u_1+h\varphi$ and $u_n=0$. The value for $u_0$ is obtained from implementing the forward-difference formula.

\section{Neural Network Method} \label{sec4}

Neural networks form the foundation of deep learning, a highly active research field in the machine learning and artificial intelligence communities \cite{gareth2013introduction}. Various architectures have gained popularity for specific applications, such as multilayer neural networks used for regression \cite{white1990connectionist}, convolutional neural networks (CNNs) applied to image classification \cite{rawat2017deep}, and recurrent neural networks (RNNs) employed for time series prediction \cite{hewamalage2021recurrent}.

In this research, we use multilayer neural networks which are commonly used for approximating continuous functions 
\cite{dissanayake1994neural,hornik1989multilayer,raissi2019physics,van2022optimally}. This section explains three neural networks that we propose to solve PDEs in different stages.

\subsection{Formulation of Neural Network 1 (NN1)} \label{formulation_NN1}

Consider a general PDE with a known parameter $\mu$ for the steady-state solution $u(\textbf{x}):\mathbb{R}^d\rightarrow\mathbb{R}$ on the domain $\Omega\in\mathbb{R}^d$ given by
\begin{align} \label{general}
\begin{split}
    N(\textbf{x},u,\mu) & = f(\textbf{x},\mu), \qquad \textbf{x} \in \Omega, \\ 
    B(\textbf{x},u,\mu) & = g(\textbf{x},\mu), \qquad \textbf{x} \in \partial\Omega,
\end{split}
\end{align}
where $N$ and $B$ are the differential operators in the interior and on the boundary. For easy calculation, we can rewrite the PDE in Eq.\ \eqref{general} into the following form
\begin{align} \label{PDE}
\begin{split}
    \mathcal{N}(\textbf{x},u,\mu) = N(\textbf{x},u,\mu) - f(\textbf{x},\mu) & = 0, \qquad \textbf{x} \in \Omega, \\ 
    \mathcal{B}(\textbf{x},u,\mu) = B(\textbf{x},u,\mu) - g(\textbf{x},\mu) & = 0, \qquad \textbf{x} \in \partial\Omega.
\end{split}
\end{align}
In case of the Bratu equation, we have $\mu=C$
while in case of the Burgers equation we have $\mu=\rho$ or $\mu=\varphi$.

In order to solve the PDE in Eq.\ \eqref{PDE} for a fixed value of the parameter $\mu$, a multilayer neural network (NN1) can be used to provide an approximate solution $u(\textbf{x},\theta)$ where $\theta$ represents all the weights of the network. If we choose $n_I$ collocation points in the interior, $\textbf{x}_i^I\in\Omega$, $1\le i \le n_I$, and $n_B$ points on the boundary, $\textbf{x}_i^B\in\partial\Omega$, $1\le i \le n_B$, then the weights can be obtained by solving the following system of nonlinear equations
\begin{align}
\begin{split}
    \mathcal{N}(\textbf{x}_i^I,u(\textbf{x}_i^I,\theta),\mu) & = 0, \qquad 1\le i \le n_I, \\
    \mathcal{B}(\textbf{x}_i^B,u(\textbf{x}_i^B,\theta),\mu) & = 0, \qquad 1\le i \le n_B,
\end{split}
\end{align}
The nonlinear system then can be solved by Newton's method \cite{chong2013introduction} or its variants. The nonlinear system also can be transformed into a loss function
\begin{multline}\label{loss1}
L = \;
\bracketround{ \frac{1}{n_I} \sum_{i=1}^{n_I} (\mathcal{N}(\textbf{x}_i^I,u(\textbf{x}_i^I,\theta),\mu))^2 } \\
+ \alpha \bracketround{ \frac{1}{n_B} \sum_{i=1}^{n_B} (\mathcal{B}(\textbf{x}_i^B,u(\textbf{x}_i^B,\theta),\mu))^2 },
\end{multline}
where the first term in the loss function represents the error in the interior and the second term represents the error on the boundary for the PDE in Eq.\ \eqref{PDE}.
In addition, a new weight $\alpha$ can be used to get a more precise solution.
The effect of choosing the value of $\alpha$ has been explained in \cite{van2022optimally}.
The loss function in Eq.\ \eqref{loss1} can be solved with SGD-type (stochastic gradient descent) algorithms such as the Adam optimizer \cite{kingma2014adam}. Other researchers use the limited-memory Broyden–Fletcher–Goldfarb–Shanno (L-BFGS) method \cite{liu1989limited}, which is a full-batch gradient-based optimization algorithm, and they showed that it also works well in solving PDEs \cite{raissi2019physics, van2022optimally}.

\begin{figure}[t]
    \centering
    \includegraphics[height=5cm]{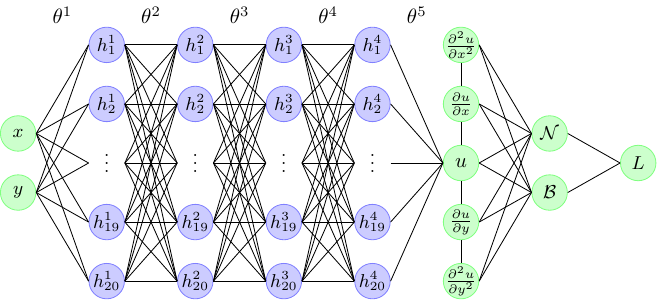}
    \caption{An example of the neural network architecture for solving PDEs.}
    \label{architectureNN1}
\end{figure}

Researchers suggest incorporating the hyperbolic tangent \cite{raissi2019physics, van2022optimally} as an activation function defined by
\begin{equation}
    \tanh(x) = \frac{e^{x}-e^{-x}}{e^{x}+e^{-x}},
\end{equation}
for the general purpose of approximating functions. Still, other differentiable activation functions also can be considered, such as the sigmoid function,
\begin{equation}
    \text{sigmoid}(x)=\frac{1}{1+e^{-x}},
\end{equation}
or the Gaussian function,
\begin{equation}
    \text{gaussian}(x)=e^{-x^2}.
\end{equation}
Note that, since the hyperbolic tangent, sigmoid, and Gaussian are in $C^{\infty}$, then for each fixed $\theta$, the map $\textbf{x} \to u(\textbf{x},\theta)$ is also in $C^{\infty}$.

In case of a neural network with two hidden layers, if we express $\theta$ as a combination of some weights and biases, $W_1$, $b_1$, $W_2$, $b_2$, $W_3$, $b_3$, then the output $u(\textbf{x},\theta)$ can be expressed by
\begin{equation} \label{NN_function}
    u(\textbf{x},\theta) = \sigma(\sigma( \textbf{x} W_1 + b_1)W_2 + b_2)W_3 + b_3,
\end{equation}
where $\sigma$ represents the activation function.
After obtaining $u(\textbf{x},\theta)$, numerical differentiation or automatic differentiation \cite{baydin2018automatic} can be used to find the derivatives that appear in Eq.\ \eqref{PDE}.

Neural networks with $p$ inputs, $l$ hidden layers with $q$ neurons in each layer, and one output has
\begin{equation}\label{numberofweight}
    (p+1)q + (l-1)(q+1)q + (q+1)
\end{equation}
weights to be optimized simultaneously. An example of the architecture of a neural network, with two inputs, four hidden layers with $20$ neurons in each layer, and one output, for solving PDEs is shown in Fig.\ \ref{architectureNN1}.

\subsubsection{Enforcing Boundary Conditions} \label{enforcing_BC}

\begin{figure}
	\centering
	\includegraphics[height=5cm]{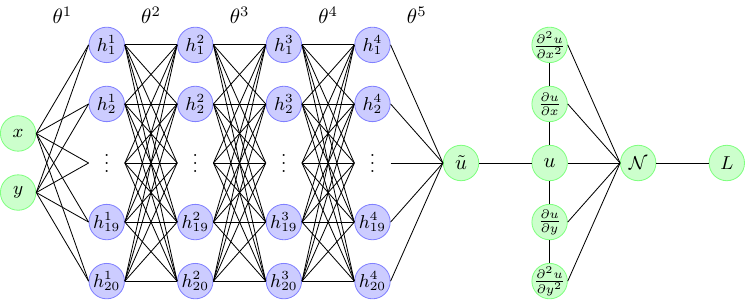}
	\caption{The architecture of neural network for solving the Bratu equation. The value of $u$ is equal to $\tilde{u} \sin(\pi x) \sin(\pi y)$.}
	\label{architectureNN3}
\end{figure}

In some cases, it is possible to enforce the output of a neural network to exactly match the boundary condition. This can be done by modifying the output of the neural network with a certain function. Using this approach, we can remove the boundary part from the nonlinear system or the loss function which makes the optimization process easier. 
We will implement this when solving the one- and two-dimensional Bratu equation. We choose to multiply the output of the neural network by $\sin(\pi x)$ for the one-dimensional case and $\sin(\pi x) \sin(\pi y)$ for the two-dimensional case. Other possible functions are $x(1-x)$ and $xy(1-x)(1-y)$, respectively. 
The illustration of the network architecture is shown in Fig.\ \ref{architectureNN3}.

\subsection{Formulation of Neural Network 2 (NN2)}
While NN1 is mainly used to solve the PDE in Eq.\ \eqref{PDE} for a fixed value of the parameter $\mu$, the multilayer neural network 2 (NN2) is its extension, to obtain a bifurcation diagram by treating $\mu$ as a variable. This requires embedding a pseudo-arclength continuation method.

\subsubsection{Pseudo-arclength Continuation} \label{subsection_arclength}
Before continuing to the next neural network, we will provide a brief explanation of pseudo-arclength continuation. This method will be used in the network later to obtain bifurcation diagrams past turning points. 

Pseudo-arclength continuation is a predictor-corrector scheme used to trace and track curves of solutions to a nonlinear system or dynamical system. It is commonly applied for studying nonlinear equations, bifurcations, and parameter sensitivity analysis \cite{allgower2003introduction,krauskopf2007numerical}.

\begin{figure}
    \centering
    \includegraphics[width=0.6\textwidth]{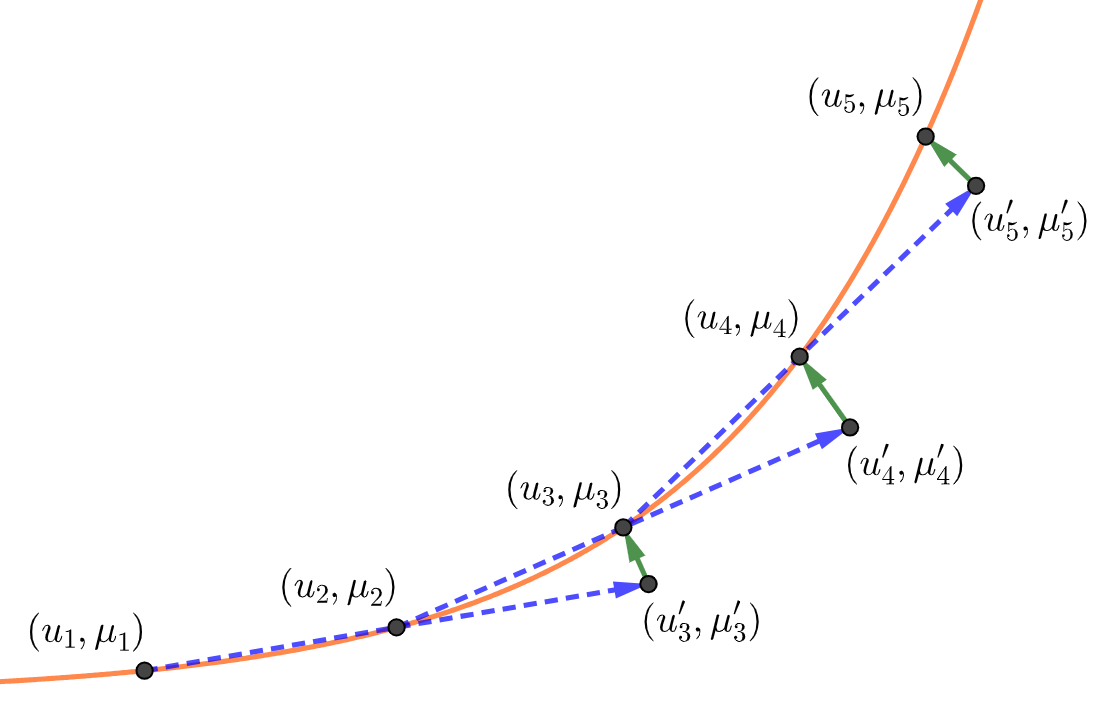}
    \caption{An example of the pseudo-arclength continuation.}
    \label{numcon}
\end{figure}

Consider a system of nonlinear equations $\mathcal{D}(u,\mu)=0$. 
We write the solution $u$ of the nonlinear system and its associated parameter $\mu$ as a pair $(u,\mu)$ for ease of notation. We can start the pseudo-arclength continuation by finding two initial solutions $(u_1,\mu_1)$ and $(u_2,\mu_2)$ of $\mathcal{D}(u,\mu)=0$ such that
\begin{equation}
	\sqrt{(\norm{u_2} - \norm{u_1})^2 + (\mu_2 - \mu_1)^2} - \delta \approx 0,
\end{equation}
where $\delta$ represents the step or length size.
The above equation has a meaning that the distance between $(\norm{u_2}, \mu_2)$ and $(\norm{u_{1}}, \mu_{1})$ must be close to $\delta$.
Then, for $k\ge3$, we can find the predicted next solution using the following predictor equation
\begin{equation}
    (u_k',\mu_k') = (u_{k-1},\mu_{k-1}) + \tau(\Dot{u}_{k-1},\Dot{\mu}_{k-1}),
\end{equation}
with
\begin{equation} \label{direction}
    (\Dot{u}_{k-1},\Dot{\mu}_{k-1}) \approx (u_{k-1},\mu_{k-1})-(u_{k-2},\mu_{k-2}),
\end{equation}
where $(\Dot{u}_{k-1},\Dot{\mu}_{k-1})$ and $\tau$ represent the direction vector and the portion of the vector, respectively \cite{allgower2003introduction,krauskopf2007numerical}.
Note that the direction vector can also be obtained using the concept of derivatives. However, computing the direction vector with Eq.\ \eqref{direction} is computationally fast and easy to implement.

Starting from the predicted solution $(u_k',\mu_k')$, a corrector method is then used to find the final solution $(u_k,\mu_k)$ satisfying
\begin{equation}
\begin{split}
    \mathcal{D}(u_k,\mu_k)=0, \\
    \sqrt{(\norm{u_k} - \norm{u_{k-1}})^2 + (\mu_k - \mu_{k-1})^2} - \delta = 0.
\end{split}
\end{equation}
A sketch of several steps of the pseudo-arclength continuation along a curve is shown in Fig.\ \ref{numcon}. Newton's method and its variants are commonly used as corrector methods \cite{allgower2003introduction,krauskopf2007numerical}.

\subsubsection{Neural Network 2 (NN2)} \label{subsubsection_NN2}

We consider the PDE in Eq.\ \eqref{PDE} as a nonlinear system which will be solved with NN2 and the pseudo-arclength continuation.
In the previous subsubsection, we have explained the principle of pseudo-arclength continuation to obtain $(u_k,\mu_k)$ knowing $(u_{k-1},\mu_{k-1})$ and $(u_{k-2},\mu_{k-2})$. 
However, because in this research we use multilayer neural networks to obtain $u(\textbf{x},\theta)$, we will use pseudo-arclength continuation to obtain a sequence of $(\theta_k,\mu_k)$ instead of $(u_k,\mu_k)$.

To implement NN2 and the pseudo-arclength continuation, we need two initial weights $(\theta_1,\mu_1)$ and $(\theta_2,\mu_2)$ where $\mu_1$, $\mu_2$ are fixed and $(u_1(\textbf{x},\theta_1),\mu_1)$, $(u_2(\textbf{x},\theta_2),\mu_2)$ are approximate solutions of the PDE in Eq.\ \eqref{PDE}. Since $\mu_1$ is fixed, we can use NN1 to obtain $(\theta_1,\mu_1)$. To obtain $(\theta_2,\mu_2)$ and maintain the pattern of the weight, we use $\theta_1$ as the initial weight for $\theta_2$ and then optimize it using NN1 to get the final weight $(\theta_2,\mu_2)$. The value of $\mu_2$ is chosen such that we have
\begin{equation}
	\sqrt{(\norm{u_2} - \norm{u_1})^2 + (\mu_2 - \mu_1)^2} - \delta \approx 0.
\end{equation}

\begin{figure}[t]
    \centering
    \includegraphics[height=5cm]{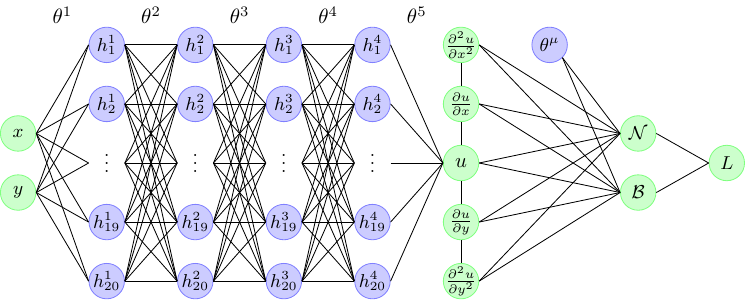}
    \caption{The architecture of NN2.}
    \label{architectureNN2}
\end{figure}

Note that in a pseudo-arclength continuation, a corrector method is used to update $\mu_k'$ into $\mu_k$, which makes the value not constant and gradually changed by the corrector method. NN1 cannot be used since it only optimizes the weight for a fixed $\mu$. Hence, we use NN2 as shown in Fig.\ \ref{architectureNN2}. Here, we treat $\mu_k\approx\theta^\mu$ as an additional parameter or weight in NN2.
Note that $\theta^\mu$ is not used to calculate the output $u$. It is only used when calculating $\mathcal{N}$ and $\mathcal{B}$.
Then, we can proceed with the pseudo-arclength continuation for $k=3,4,\dots$ by a predictor
\begin{equation}
    (\theta_k',(\theta^\mu_k)') = (\theta_{k-1},\mu_{k-1}) + \tau(\Dot{\theta}_{k-1},\Dot{\mu}_{k-1}),
\end{equation}
with
\begin{equation}
    (\Dot{\theta}_{k-1},\Dot{\mu}_{k-1}) \approx (\theta_{k-1},\mu_{k-1})-(\theta_{k-2},\mu_{k-2}).
\end{equation}
We choose $\tau=1$ for this research. Starting from the predicted weight $(\theta_k',(\theta^\mu_k)')$, we need to find the final weight $(\theta_k,\theta^\mu_k)$ satisfying the following system of nonlinear equations
\begin{align}\label{system2}
\begin{split}
    \mathcal{N}(\textbf{x}_i^I,u_k(\textbf{x}_i^I,\theta_k),\theta^\mu_k) = 0, \qquad 1\le i \le n_I, \\
    \mathcal{B}(\textbf{x}_i^B,u_k(\textbf{x}_i^B,\theta_k),\theta^\mu_k) = 0, \qquad 1\le i \le n_B, \\
    \sqrt{(\norm{u_k} - \norm{u_{k-1}})^2 + (\theta^\mu_k - \mu_{k-1})^2} - \delta = 0.
\end{split}
\end{align}
This system of equations can be transformed into an equivalent loss function
\begin{multline}
L =
\bracketround{ \frac{1}{n_I} \sum_{i=1}^{n_I} (\mathcal{N}(\textbf{x}_i^I,u_k(\textbf{x}_i^I,\theta_k),\theta^\mu_k))^2 } \\
+ \alpha \bracketround{ \frac{1}{n_B} \sum_{i=1}^{n_B} (\mathcal{B}(\textbf{x}_i^B,u_k(\textbf{x}_i^B,\theta_k),\theta^\mu_k))^2 } \\
+ \beta \bracketround{ \sqrt{(\norm{u_k} - \norm{u_{k-1}})^2 + (\theta^\mu_k - \mu_{k-1})^2} - \delta }^2,
\end{multline}
where a new weight $\beta$ can be added to get a more precise solution.

After we find the final weight $\theta_k\cup\bracketcurly{\theta^\mu_k}$, we set $\mu_k=\theta^\mu_k$ and hence we obtain $(\theta_k,\mu_k)$. The process is repeated until we cover all possible values of $\mu$ and its associated solutions. After we obtain the sequence of $(\theta_k,\mu_k)$, we can compute $(u_k,\mu_k)$ to get the bifurcation diagram.

Using NN2, we can find a sequence of solutions without being constrained by the value of the turning point $\mu^*$. Additionally, it can determine the correct direction or path of subsequent solutions if the step size $\delta$ is sufficiently small \cite{krauskopf2007numerical}. Moreover, NN2 can significantly reduce computational time since it eliminates the need to determine optimal weights for the next solutions from scratch.
It is worth noting that we can also use parameter continuation, where the nonlinear system is solved by iteratively adjusting the fixed value of $\mu$ manually for several steps to obtain more consistent weights before using pseudo-arclength continuation. Additionally, note that the proposed concept bears some resemblance to the so-called evolutionary neural networks (ENN) \cite{du2021evolutional}, where its weights evolve over time.

\subsection{Formulation of Neural Network 3 (NN3)}

Consider a general eigenvalue problem with an unknown parameter $\lambda$ for the function $v(\textbf{x}):\mathbb{R}^d\rightarrow\mathbb{R}$ on the domain $\Omega\in\mathbb{R}^d$ given by
\begin{align} \label{eigenvalueproblem_NN}
\begin{split}
    \mathcal{N}(\textbf{x},v,\mu,\lambda) = N(\textbf{x},v,\mu) - \lambda v & = 0, \qquad \textbf{x} \in \Omega, \\ 
    \mathcal{B}(\textbf{x},v,\mu) = B(\textbf{x},v,\mu) - g(x,\mu) & = 0, \qquad \textbf{x} \in \partial\Omega,
\end{split}
\end{align}
where $\mu$ is a fixed parameter. In this case, $\lambda$ represents an eigenvalue and $v$ represents its associated eigenfunction.

Using the Rayleigh quotient (see Eq.\ \eqref{Rayleighquotient} or \eqref{rayleigh_burgers}), we can approximate $\lambda$ as
\begin{equation}
    \lambda = \frac{\innerproduct{ N(\textbf{x},v,\mu), v }}{\innerproduct{ v, v }},
\end{equation}
and therefore Eq.\ \eqref{eigenvalueproblem_NN} can be transformed into
\begin{align} \label{eigenvalueproblem2}
\begin{split}
    \mathcal{N}_2(\textbf{x},v,\mu) & = 0, \qquad \textbf{x} \in \Omega, \\ 
    \mathcal{B}_2(\textbf{x},v,\mu) & = 0, \qquad \textbf{x} \in \partial\Omega.
\end{split}
\end{align}
The main advantage of using the Rayleigh quotient is reducing the complexity of the eigenvalue problem.

It is clear that the problem in Eq.\ \eqref{eigenvalueproblem2} is similar to the general PDE in Eq.\ \eqref{PDE}. Hence, we can use a multilayer neural network (NN3) which has a similar architecture as NN1 to solve this problem. In addition, to solve the eigenvalue problem, we set $\norm{v}_2=1$ as an additional condition of the eigenfunction in NN3. Then the solution can be obtained by solving the following system of nonlinear equations 
\begin{align}
\begin{split}
    \mathcal{N}_2(\textbf{x}_i^I,v(\textbf{x}_i^I,\theta),C) & = 0, \qquad 1\le i \le n_I, \\
    \mathcal{B}_2(\textbf{x}_i^B,v(\textbf{x}_i^B,\theta),C) & = 0, \qquad 1\le i \le n_B, \\
    \norm{v}_2 - 1 & = 0.
\end{split}
\end{align}
After transforming the nonlinear system into the associated loss function, we will get
\begin{multline}\label{loss2}
L =
\bracketround{ \frac{1}{n_I} \sum_{i=1}^{n_I} (\mathcal{N}_2(\textbf{x}_i^I,v(\textbf{x}_i^I,\theta),C))^2 } \\
+ \alpha \bracketround{ \frac{1}{n_B} \sum_{i=1}^{n_B} (\mathcal{B}_2(\textbf{x}_i^B,v(\textbf{x}_i^B,\theta),C))^2 }
+ \gamma \bracketround{ \norm{v}_2 - 1 }^2,
\end{multline}
where a new weight $\gamma$ can be added to the loss function to get a more precise solution.

Although, in the implementation, we will use the above concept, note that we can also use a multilayer neural network similar to NN2 to solve the eigenvalue problem in Eq.\ \eqref{eigenvalueproblem_NN} directly. This is done by treating $\lambda\approx\theta^{\lambda}$ as a parameter or weight in the network. 

Since the purpose of this part is to obtain the largest eigenvalue, following the Sturm-Liouville theory \cite{trench2013elementary}, we seek an eigenfunction without any zeros in the domain (i.e., that does not intersect the $x$ axis or the $xy$ plane for the one- and two-dimensional cases). Using the largest eigenvalue, we can determine whether the solution is linearly stable or unstable.

\section{Experimental Results and Discussions} \label{sec5}

In this section, we will compare the proposed neural network (NN), explained in the previous section, with the finite difference method (FD). For both approaches, the Levenberg-Marquardt method \cite{marquardt1963algorithm}, a variant of Newton's method, is used to solve the nonlinear systems arising from the equation and as a corrector method for constructing bifurcation diagrams. The Levenberg-Marquardt method is available in the 'fsolve' function in MATLAB.
To obtain an approximate solution in this research, we use a neural network with only two hidden layers. Each layer consists of 5 neurons for the one-dimensional Bratu equation and the one-dimensional Burgers equation, or 10 neurons for the two-dimensional Bratu equation. Using Eq.\ \eqref{numberofweight}, the number of variables (weights) of NN is 46 for the one-dimensional case and 151 for the two-dimensional case. These will be considered as the variables that need to be optimized by Newton's method.
Note that when we refer to NN with $n=100$, it means that NN will use the same collocation or grid points as used by FD with $n=100$, in both the one- and two-dimensional problems.
To carry out the experiments, we used a laptop with an Intel Core i7-1185G7 processor. The experiments were implemented in MATLAB R2022b.

\subsection{One-Dimensional Bratu Equation}

For solving the one-dimensional Bratu equation with FD, we use $n=100$, where the domain is divided into 100 equal subintervals. Thus, FD needs to optimize 99 variables, corresponding to the 99 values of the solution at the interior points. 
To compare the methods fairly, we use the same collocation or grid points for NN, namely $n=100$.  The Gaussian function is chosen as the activation function for the network. 
When implementing this scheme, the size of the Jacobian matrix for FD is $99\times99$, while the Jacobian matrix for NN is $46\times99$.

First, we present the results for solving the one-dimensional Bratu equation to obtain the lower solution for $C=3.513$ and the upper solution for $C=1$. Since the analytic solution is available, we will use it as a reference. Solving Eq.\ \eqref{coshtheta} for $C=3.513$
and $C=1$ results in $\omega$ values of 1.184367 and 2.734676, respectively.
For FD, we set $u_i=0$, for $i=1,\dots,n-1$ as the initial values to obtain the lower solution for $C=3.513$. To obtain the upper solution for $C=1$, the initial values are obtained from $u_i=4\sin{\pi x_i}$. 
For NN, we set all the weights randomly within the range $(-0.01,0.01)$ to obtain the lower solution. To obtain the upper solution, we set the bias $b_3=4$.

\begin{figure}
	\centering
	\begin{subfigure}[b]{0.49\textwidth}
		\centering
		\includegraphics[width=\textwidth]{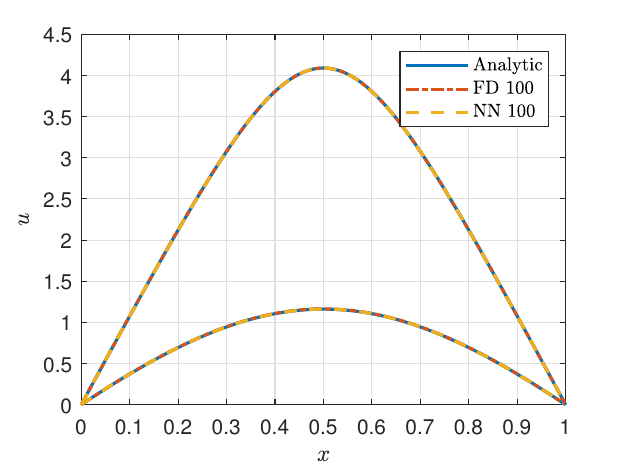}
		\caption{}
		\label{fig_bratu_solution1}
	\end{subfigure}
	\begin{subfigure}[b]{0.49\textwidth}
		\centering
		\includegraphics[width=\textwidth]{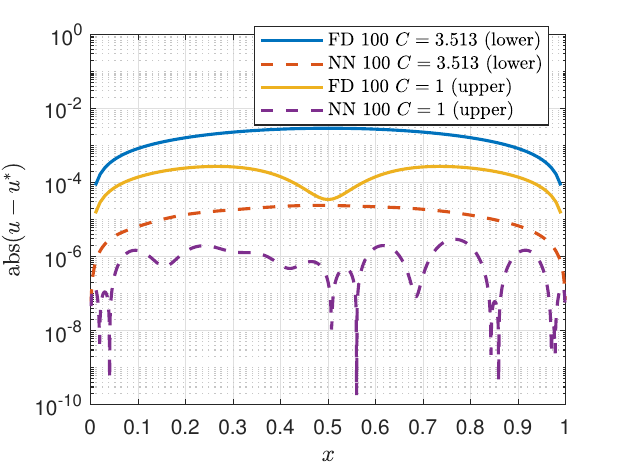}
		\caption{}
		\label{fig_bratu_solution2}
	\end{subfigure}
	\caption{(a) The lower solutions for $C=3.513$ and the upper solution for $C=1$, (b) the absolute error of solutions from NN and FD compared to the analytic solution.}
	\label{fig_bratu_solution}
\end{figure}

\begin{table}[H]
    \centering
    \caption{Mean squared error (MSE) of solutions obtained from NN and FD compared to the analytic solution.}
    \label{tab:1dsolution}
    \begin{tabular}{|c|c|c|c|} 
        \hline
        Lower/Upper Solution & $C$ & FD & NN \\ 
        \hline
        Lower Solution & 1.0 & 2.6211e-12 & \textcolor{blue}{1.3717e-13} \\ 
        Lower Solution & 1.5 & 1.4649e-11 & \textcolor{blue}{2.5441e-15} \\
        Lower Solution & 2.0 & 2.0818e-11 & \textcolor{blue}{1.2807e-12} \\
        Lower Solution & 2.5 & 1.5417e-10 & \textcolor{blue}{1.5873e-14} \\
        Lower Solution & 3.0 & 8.7843e-10 & \textcolor{blue}{1.4479e-13} \\
        \hline
        Upper Solution & 1.0 & 4.3531e-08 & \textcolor{blue}{1.4403e-12} \\
        Upper Solution & 1.5 & 1.4874e-08 & \textcolor{blue}{1.2240e-12} \\
        Upper Solution & 2.0 & 1.7309e-08 & \textcolor{blue}{9.7839e-13} \\
        Upper Solution & 2.5 & 9.6474e-09 & \textcolor{blue}{5.4528e-13} \\
        Upper Solution & 3.0 & 7.9989e-09 & \textcolor{blue}{9.4078e-12} \\
        \hline
    \end{tabular}
\end{table}

The solutions obtained from NN and FD are shown in Fig.\ \ref{fig_bratu_solution1}, and the absolute errors compared to the analytic solution are given in Fig.\ \ref{fig_bratu_solution2}. From the figures, we observe that NN provides more accurate solutions compared to FD for both cases, $C=3.513$ and $C=1$. 
Additionally, in Tab.\ \ref{tab:1dsolution}, the accuracy of NN and FD is tested for several values of $C$. Analytic solutions are used as references to compute the mean squared error (MSE), with the results showing smaller errors highlighted in blue. From the table, we see that NN outperforms FD in all 10 cases.
In terms of computational time, FD takes around 0.1 seconds, while NN takes around 0.5 seconds to complete the search process.

\begin{figure}
    \centering
    \begin{subfigure}[b]{0.49\textwidth}
        \centering
        \includegraphics[width=\textwidth]{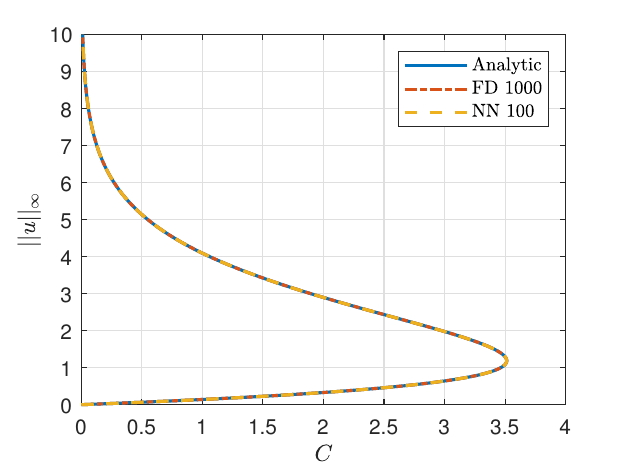}
        \caption{}
        \label{fig_bratu_bifurcation1}
    \end{subfigure}
    \begin{subfigure}[b]{0.49\textwidth}
        \centering
        \includegraphics[width=\textwidth]{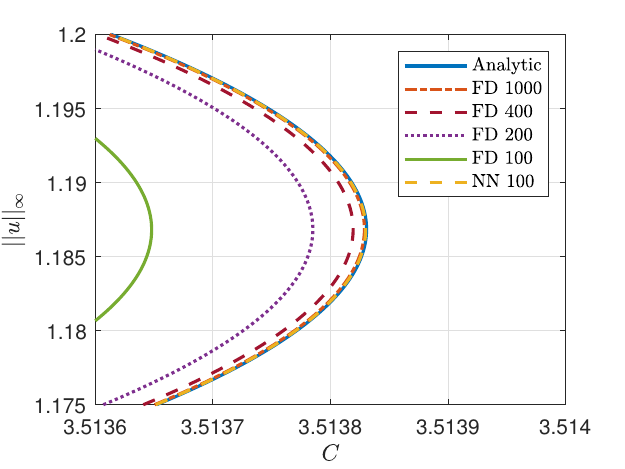}
        \caption{}
        \label{fig_bratu_bifurcation2}
    \end{subfigure}
    \caption{(a) Bifurcation diagrams obtained by FD with $n=1000$ and NN with $n=100$, (b) zoom of the bifurcation diagrams close to the turning point.}
    \label{fig_bratu_bifurcation}
\end{figure}

Next, we compare bifurcation diagrams constructed using NN, FD, and the analytic solution. Fig.\ \ref{fig_bratu_bifurcation1} shows the full bifurcation diagrams obtained from NN with $n=100$ and FD with $n=1000$. Both are similar to the analytic solution, with barely visible differences.
To assess the accuracy more clearly, we zoom in on the results around the turning point, as shown in Fig.\ \ref{fig_bratu_bifurcation2}. In this case, we consider $n=100$, $200$, $400$, and $1000$ for FD. It is clear from the figure that the turning points obtained from FD become less accurate as the number of grid points decreases. On the other hand, the bifurcation diagram from NN with $n=100$ is more accurate compared to those obtained from FD, even more accurate than FD with $n=1000$.

\begin{figure}[t]
    \centering
    \begin{subfigure}[b]{0.49\textwidth}
        \centering
        \includegraphics[width=\textwidth]{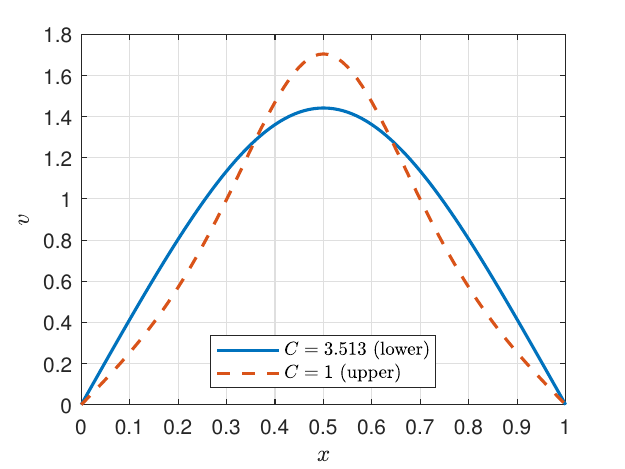}
        \caption{}
        \label{fig_bratu_stability1}
    \end{subfigure}
    \begin{subfigure}[b]{0.49\textwidth}
        \centering
        \includegraphics[width=\textwidth]{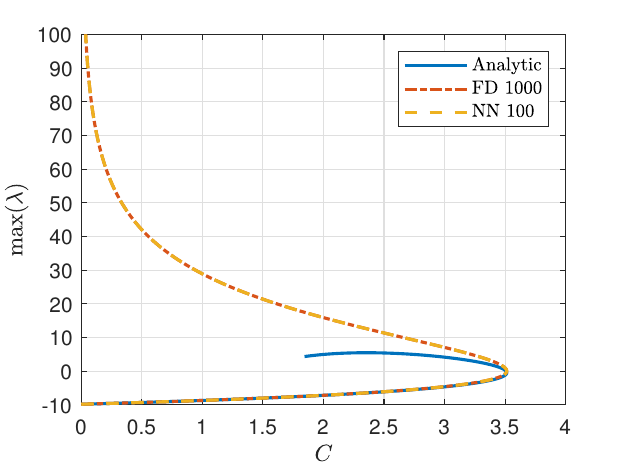}
        \caption{}
        \label{fig_bratu_stability2}
    \end{subfigure}
    
    \caption{(a) Eigenfunctions of the lower solution for $C=3.513$ and the upper solution for $C=1$ obtained from NN with $n=100$, (b) largest eigenvalues against $C$ for the one-dimensional Bratu equation.}
    \label{fig_bratu_stability}
\end{figure}

Next, we discuss the use of NN for solving the eigenvalue problem in Eq.\ \eqref{eigenvalueproblem} for the one-dimensional case. NN is used to solve the eigenvalue problem for the upper solution when $C=1$, which gives $\max(\lambda)\approx28.872$. This result is similar to the largest eigenvalue obtained from the Jacobian matrix of FD, $\max(\lambda)\approx28.873$.
In Fig.\ \ref{fig_bratu_stability1}, the computed eigenfunctions from NN for the lower solution when $C=3.513$ and the upper solution when $C=1$ are shown.

In Fig.\ \ref{fig_bratu_stability2}, the complete graphs of the largest eigenvalue, $\max(\lambda)$, against $C$ from NN with $n=100$ and FD with $n=1000$ are shown. From the figure, we observe that NN and FD provide very similar results as the lines overlap each other. This indicates that the neural network proposed in this research can solve the linear stability problems effectively.
Furthermore, \(\max(\lambda)\) at the turning point is very close to 0. This informs us that \(\max(\lambda)\) for all lower solutions is negative, and \(\max(\lambda)\) for all upper solutions is positive. Hence, it can be concluded that all lower solutions are linearly stable, and all upper solutions are linearly unstable.
In addition, we provide an asymptotic analysis in \ref{appendix1} (i.e., Eq.\ \eqref{app_exp}) to approximate the largest eigenvalue, as shown in Fig.\ \ref{fig_bratu_stability2}. Note that even though we only used the first three terms of the asymptotic expansion, our analysis is able to capture the turning point rather accurately.

\subsection{Two-Dimensional Bratu Equation}

\begin{figure}
    \centering
    \begin{subfigure}[b]{0.49\textwidth}
        \centering
        \includegraphics[width=\textwidth]{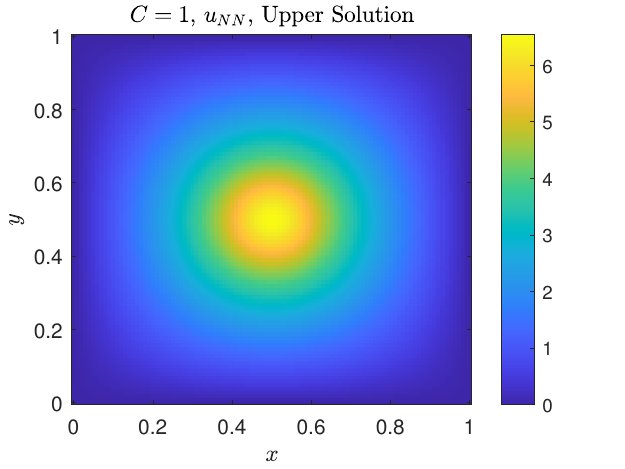}
        \caption{}
        \label{fig_bratu2_solution1}
    \end{subfigure}
    \begin{subfigure}[b]{0.49\textwidth}
        \centering
        \includegraphics[width=\textwidth]{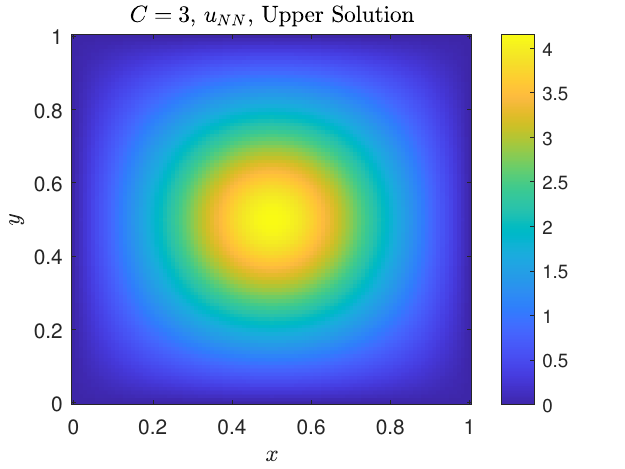}
        \caption{}
        \label{fig_bratu2_solution4}
    \end{subfigure}
    
	\;
    
    \begin{subfigure}[b]{0.49\textwidth}
        \centering
        \includegraphics[width=\textwidth]{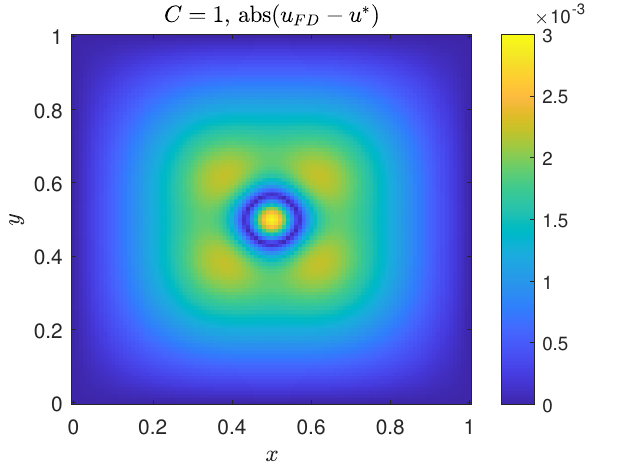}
        \caption{}
        \label{fig_bratu2_solution2}
    \end{subfigure}
    \begin{subfigure}[b]{0.49\textwidth}
        \centering
        \includegraphics[width=\textwidth]{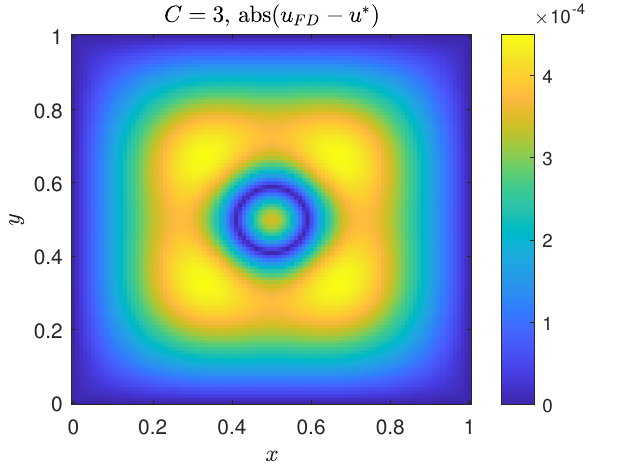}
        \caption{}
        \label{fig_bratu2_solution5}
    \end{subfigure}
    
    \;
    
    \begin{subfigure}[b]{0.49\textwidth}
        \centering
        \includegraphics[width=\textwidth]{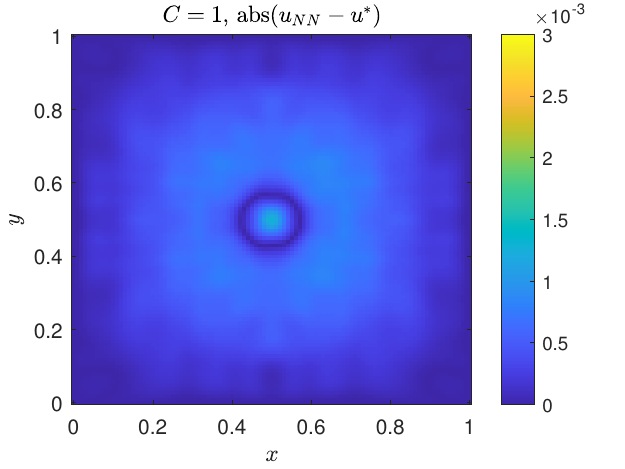}
        \caption{}
        \label{fig_bratu2_solution3}
    \end{subfigure}
    \begin{subfigure}[b]{0.49\textwidth}
        \centering
        \includegraphics[width=\textwidth]{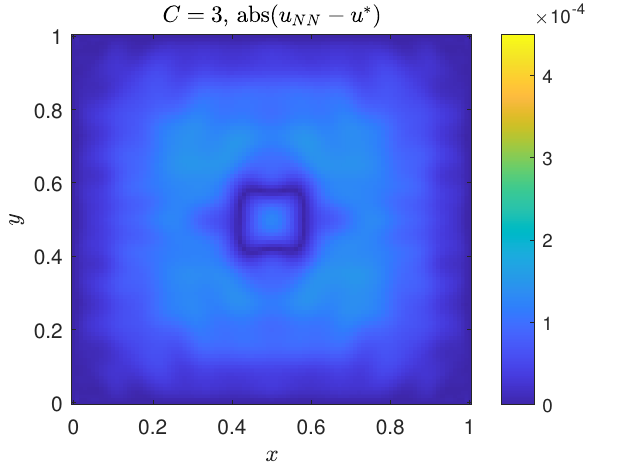}
        \caption{}
        \label{fig_bratu2_solution6}
    \end{subfigure}
    
    \caption{(a) Upper solutions from NN with $n=100$ for $C=1$ and (b) $C=3$, 
    	(c) the absolute error from FD for $C=1$ and (d) $C=3$,
    	(e) the absolute error from NN for $C=1$ and (f) $C=3$.}
    \label{fig_bratu2_solution}
\end{figure}

For solving the two-dimensional Bratu equation with FD, we consider $n=100$, which means both intervals on the $x$-axis and $y$-axis are divided into 100 equal subintervals. With these grid points, FD needs to optimize $99\times 99=9801$ variables. To compare the methods fairly, we use the same collocation or grid points for NN, which we will call NN with $n=100$. NN with two hidden layers and 10 neurons in each layer will have 151 weights according to Eq.\ \eqref{numberofweight}. The Gaussian function is chosen as the activation function for the network.
Furthermore, the sizes of the Jacobian matrices for FD and NN will be $99^2\times99^2$ and $151\times99^2$, respectively. In this case, FD will need around 65 times more random access memory (RAM) than NN.

\begin{figure}
    \centering
    \begin{subfigure}[b]{
    		0.49\textwidth}
        \centering
        \includegraphics[width=\textwidth]{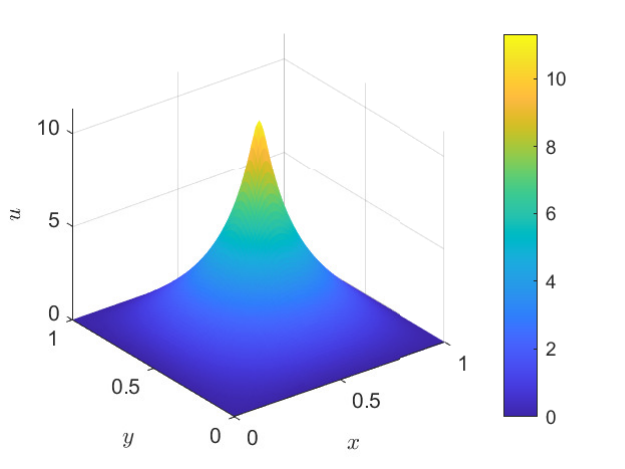}
        \caption{}
        \label{Fig_2d_sharp}
    \end{subfigure}
    \begin{subfigure}[b]{0.49\textwidth}
        \centering
        \includegraphics[width=\textwidth]{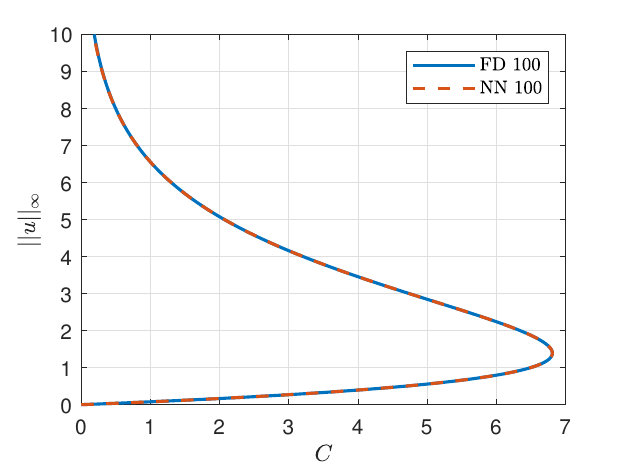}
        \caption{}
        \label{Fig_2d_bifurcation1}
    \end{subfigure}
    
    \begin{subfigure}[b]{0.49\textwidth}
        \centering
        \includegraphics[width=\textwidth]{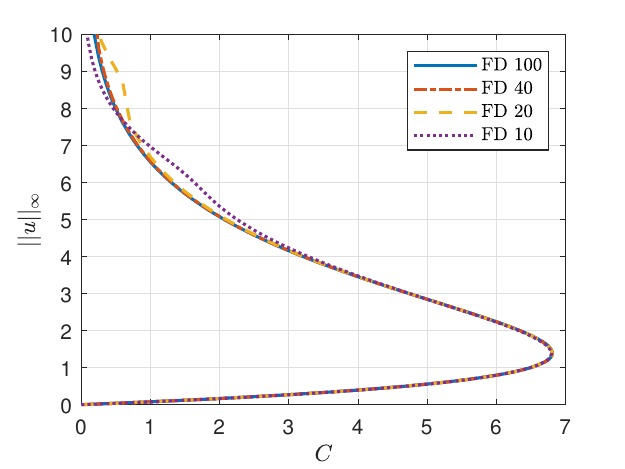}
        \caption{}
        \label{Fig_2d_bifurcation_FD}
    \end{subfigure}
    \begin{subfigure}[b]{0.49\textwidth}
        \centering
        \includegraphics[width=\textwidth]{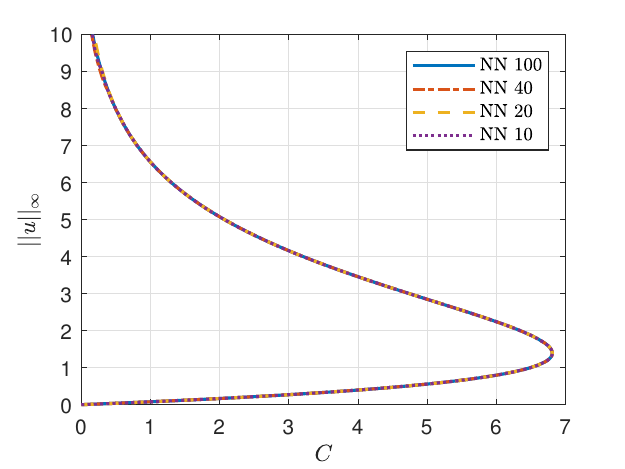}
        \caption{}
        \label{Fig_2d_bifurcation_NN}
    \end{subfigure}
    \caption{(a) Upper solutions for the two-dimensional Bratu equation have a sharp gradient at $(x,y)=(0.5,0.5)$. Shown is an example for $C=0.1$, (b) bifurcation diagrams from NN and FD with $n=100$, (c) bifurcation diagrams obtained from FD with different grid points, (d) bifurcation diagrams obtained from NN with different grid points.}
    \label{Fig_2d_bifurcation}
\end{figure}

First, we illustrate the obtained upper solutions from NN with $n=100$ for $C=1$ and $C=3$ in Figs.\ \ref{fig_bratu2_solution1} and \ref{fig_bratu2_solution4}. To compare the results of NN and FD, both with $n=100$, we consider the results from FD with $n=200$ as the true solutions since the analytic solution for the two-dimensional Bratu equation is not available. In Figs.\ \ref{fig_bratu2_solution2} and \ref{fig_bratu2_solution5}, we show the absolute error between FD and the true solutions, while in Figs.\ \ref{fig_bratu2_solution3} and \ref{fig_bratu2_solution6}, the absolute error from NN is shown. 
These two examples show that NN gives a smaller error compared to FD with the same grid points. Regarding computational time, FD takes around 40 seconds, while NN takes around 150 seconds to complete the search process. Moreover, as shown in Fig.\ \ref{Fig_2d_sharp}, upper solutions for the two-dimensional Bratu equation have sharp gradients. This shape will pose problems for FD if it uses a small number of grid points.

Next, the comparison of bifurcation diagrams from NN and FD, both with $n=100$, is shown in Fig.\ \ref{Fig_2d_bifurcation1}. Additionally, Figs.\ \ref{Fig_2d_bifurcation_FD} and \ref{Fig_2d_bifurcation_NN} consider different numbers of grid points: $n=10,20,40,100$.
Both FD and NN result in similar bifurcation diagrams for $n=100$, as indicated in Fig.\ \ref{Fig_2d_bifurcation1}. As previously mentioned, since the upper solutions have sharp gradients, FD gives inaccurate bifurcation diagrams for smaller numbers of grid points, i.e. $n=20$ and $n=10$. On the other hand, NN produces very good bifurcation diagrams, even when a small number of grid points is used.

\begin{figure}
    \centering
    \includegraphics[width=0.5\textwidth]{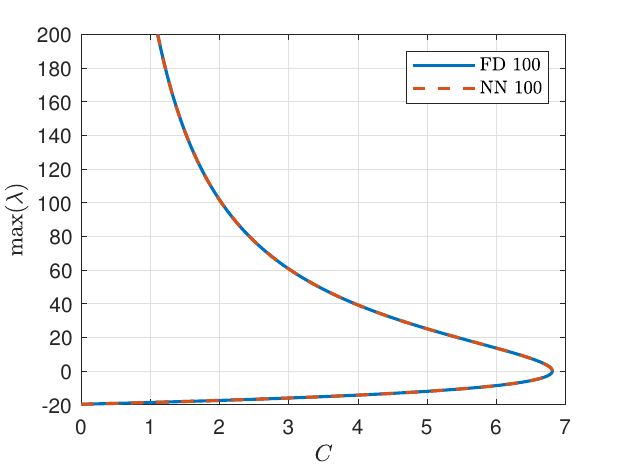}
    \caption{The largest eigenvalue against $C$ for the two-dimensional Bratu equation.}
    \label{2d_stabilty}
\end{figure}

Last but not least, the graphs of the largest eigenvalue against \(C\) from NN and FD, both with $n=100$, are shown in Fig.\ \ref{2d_stabilty}. In this case, NN and FD also provide very similar results. As in the one-dimensional case, \(\max(\lambda)\approx0\) is attained at the turning point. Therefore, for the two-dimensional Bratu equation, all lower solutions are linearly stable, and all upper solutions are linearly unstable.

\subsection{One-Dimensional Burgers Equation} \label{results_burgers}

\begin{figure}[t]
	\centering
	\begin{subfigure}[b]{0.49\textwidth}
		\centering
		\includegraphics[width=\textwidth]{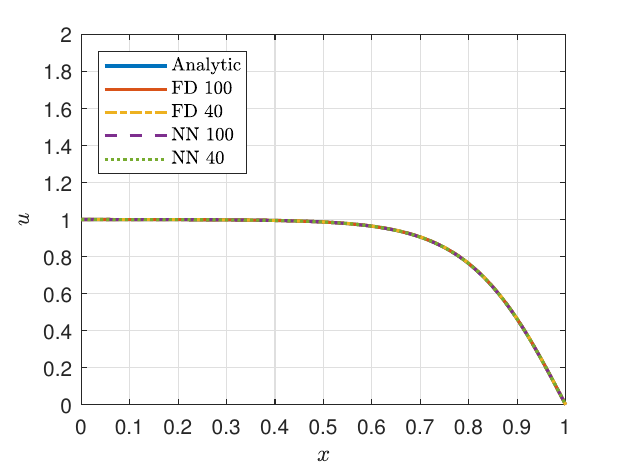}
		\caption{}
		\label{fig_burgers_solution1}
	\end{subfigure}
	\begin{subfigure}[b]{0.49\textwidth}
		\centering
		\includegraphics[width=\textwidth]{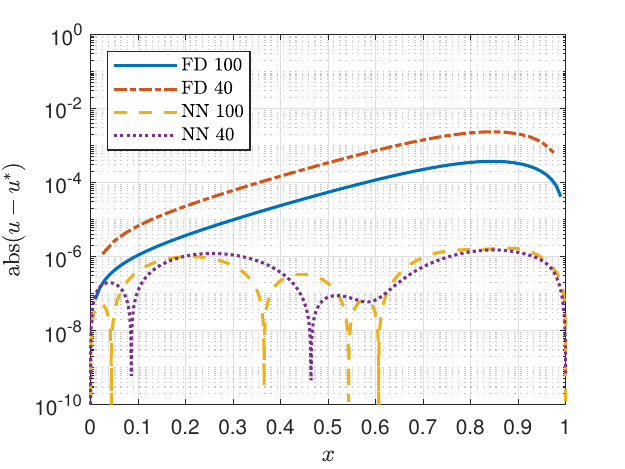}
		\caption{}
		\label{fig_burgers_solution2}
	\end{subfigure}
	
	\begin{subfigure}[b]{0.49\textwidth}
		\centering
		\includegraphics[width=\textwidth]{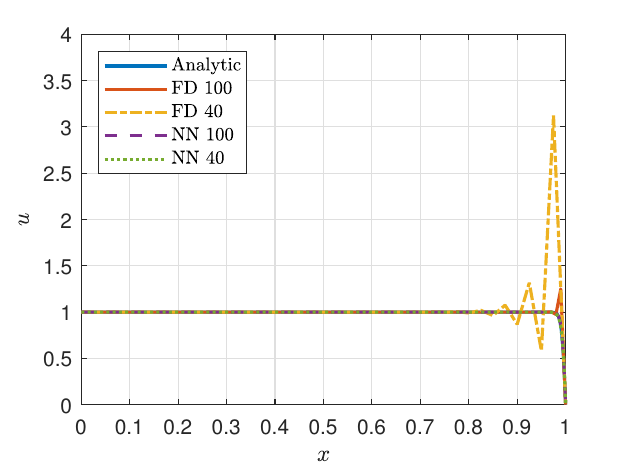}
		\caption{}
		\label{fig_burgers_solution3}
	\end{subfigure}
	\begin{subfigure}[b]{0.49\textwidth}
		\centering
		\includegraphics[width=\textwidth]{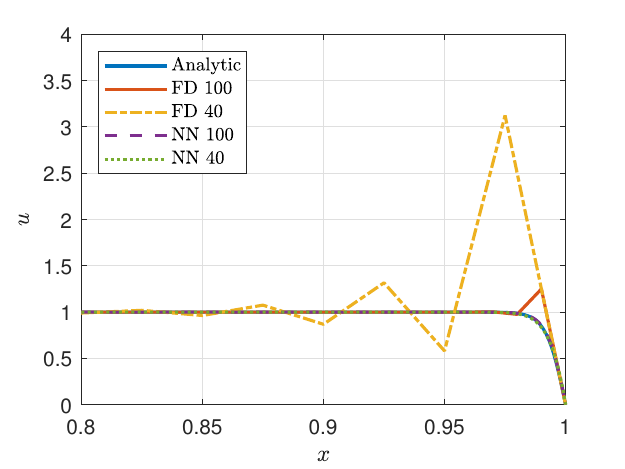}
		\caption{}
		\label{fig_burgers_solution4}
	\end{subfigure}
	\caption{(a) Solutions obtained from NN, FD, and the analytic solution for $\nu=0.1$, (b) the absolute error from NN and FD compared to the analytic solution for $\nu=0.1$, (c) solutions obtained from NN, FD, and the analytic solution for $\nu=0.004$, (d) solutions for $x\in[0.8,1]$.}
	\label{fig_burgers_solution}
\end{figure}

First, we compare the performance of NN and FD for solving the Burgers equation with the Dirichlet boundary condition in Eq.\ \eqref{burgers_BC1}. For this experiment, we choose $n=40$ and $n=100$ for both FD and NN.
We use NN with a similar setting as used in the one-dimensional Bratu equation, however we implement the tanh activation function.

We display the results for solving the one-dimensional Burgers equation for $\nu=0.1$ and $\nu=0.004$. Since the analytic solution is available, we will also use it as a reference. The value of $\rho$ is obtained from Eq.\ \eqref{burgers_rho}. For FD, we set $u_i=1$ as the initial value, while for NN, we initialize $b_3=1$ and set all other weights randomly within the range $(-0.01,0.01)$.

In Fig.\ \ref{fig_burgers_solution1}, we observe that both FD and NN produce satisfactory solutions for $\nu=0.1$. To assess accuracy more precisely, we provide the absolute error between the obtained solutions and the analytic solution in Fig.\ \ref{fig_burgers_solution2}. Generally, NN exhibits smaller errors compared to FD. Furthermore, we notice that the error in NN with $n=40$ has a similar magnitude, around $10^{-6}$, compared to NN with $n=100$.
In terms of computational time for $n=100$, FD takes approximately 0.1 seconds, while NN takes around 1.5 seconds to complete the search process.

Next, we plot the solutions for $\nu=0.004$ in Figs.\ \ref{fig_burgers_solution3} and \ref{fig_burgers_solution4}. This parameter yields a solution with sharper gradients when $x$ is close to 1. Unfortunately, FD with $n=40$ and $n=100$ fails to converge to the correct solution.
This example illustrates that in some cases where the solution exhibits sharp gradients, FD requires more points to converge accurately. Conversely, NN is capable of producing satisfactory solutions, both for $n=40$ and $n=100$.

\begin{figure}
	\centering
	\begin{subfigure}[b]{0.49\textwidth}
		\centering
		\includegraphics[width=\textwidth]{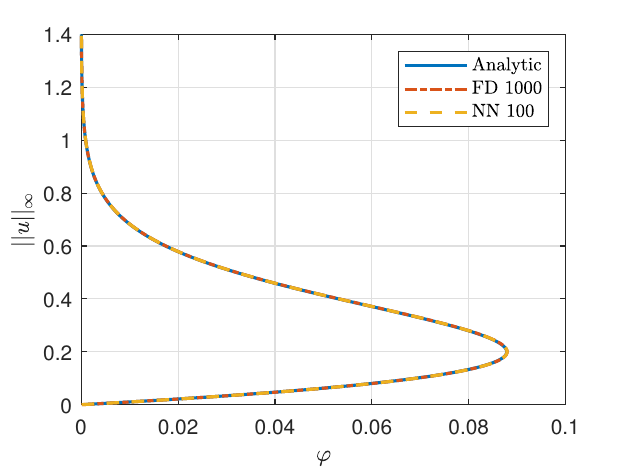}
		\caption{}
		\label{fig_burgers_bifurcation1}
	\end{subfigure}
	\begin{subfigure}[b]{0.49\textwidth}
		\centering
		\includegraphics[width=\textwidth]{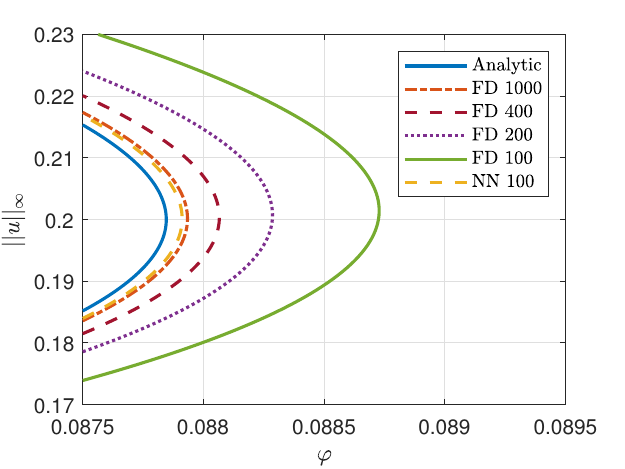}
		\caption{}
		\label{fig_burgers_bifurcation2}
	\end{subfigure}
	\caption{(a) Bifurcation diagrams from NN with $n=100$ and FD with $n=1000$ compared to the analytic solution, (b) zoom of the bifurcation diagrams close to the turning point.}
	\label{fig_burgers_bifurcation}
\end{figure}

We compare bifurcation diagrams for the Burgers equation with the mixed boundary condition in Eq.\ \eqref{burgers_BC2}, constructed by NN, FD, and the analytic solution. We set $\nu=0.1$ since, from the previous example, solutions from both NN and FD are close to the analytic solution.
First, we consider NN with $n=100$ and FD with $n=1000$ to construct the full bifurcation diagrams of the Burgers equation, as shown in Fig.\ \ref{fig_burgers_bifurcation1}. Both NN and FD produce results similar to the analytic solution.
Next, we zoom in on the bifurcation around the turning point, as shown in Fig.\ \ref{fig_burgers_bifurcation2}. In addition, FD with $n=100$, $200$, $400$, and $1000$ is used in this case. Surprisingly, the bifurcation diagram from NN with $n=100$ is slightly closer to the analytic solution than the one obtained from FD with $n=1000$. The turning point obtained from FD becomes less accurate as we decrease the number of grid points. This demonstrates that the bifurcation diagram obtained from NN is more accurate compared to FD.

\begin{figure}
	\centering
	\begin{subfigure}[b]{0.49\textwidth}
		\centering
		\includegraphics[width=\textwidth]{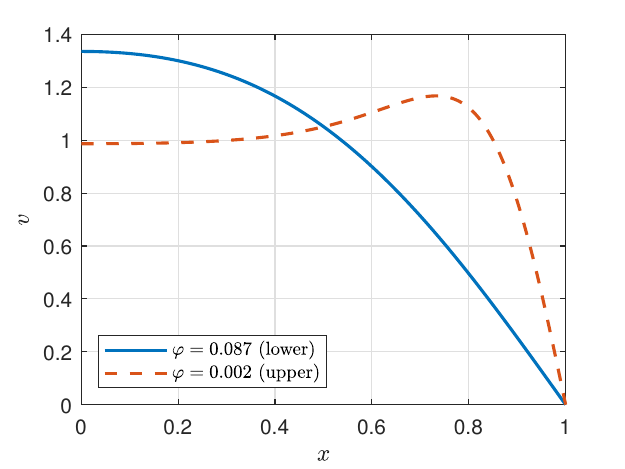}
		\caption{}
		\label{fig_burgers_stability1}
	\end{subfigure}
	\begin{subfigure}[b]{0.49\textwidth}
		\centering
		\includegraphics[width=\textwidth]{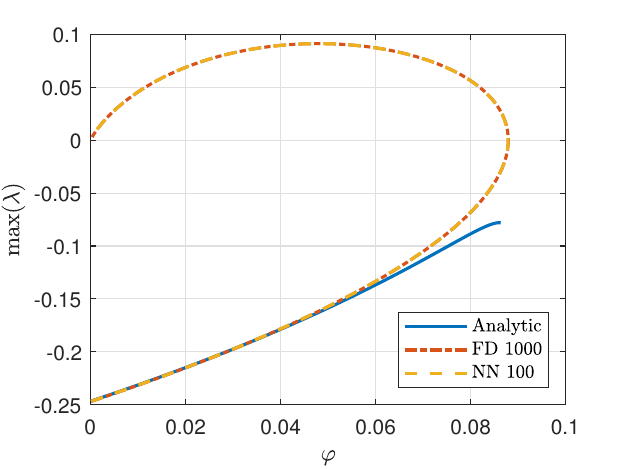}
		\caption{}
		\label{fig_burgers_stability2}
	\end{subfigure}
	\caption{(a) Eigenfunctions of the lower solution for $\varphi=0.087$ and the upper solution for $\varphi=0.002$, (b) largest eigenvalues against $\varphi$ for the one-dimensional Burgers equation.}
	\label{fig_burgers_stability}
\end{figure}

Next, we discuss the use of NN for solving the eigenvalue problem in Eq.\ \eqref{eigenvalueproblem_burgers}. For the lower solution when $\varphi=0.06$, NN yields $\max(\lambda)\approx0.08727$, which is similar to the result obtained from the largest eigenvalue of the Jacobian matrix of FD, $\max(\lambda)\approx0.08729$. 
In Fig.\ \ref{fig_burgers_stability1}, we plot the obtained eigenfunctions from NN with $n=100$ for the lower solution when $\varphi=0.087$ and the upper solution when $\varphi=0.02$. 

In Fig.\ \ref{fig_burgers_stability2}, the complete graphs of the largest eigenvalue, $\max(\lambda)$, against $\varphi$ from NN with $n=100$ and FD with $n=1000$ are provided. From the figure, it is evident that NN and FD yield very similar results, as the lines overlap each other. This indicates that the neural network proposed in this research can effectively solve linear stability problems.
Furthermore, the value of $\max(\lambda)$ at the turning point is very close to 0. This implies that $\max(\lambda)$ for all lower solutions is negative, while $\max(\lambda)$ for all upper solutions is positive. Hence, based on this result, we can conclude that all lower solutions are linearly stable, and all upper solutions are linearly unstable.
In addition, we provide an asymptotic analysis in \ref{appendix1} to approximate the largest eigenvalue (see Eq.\ \eqref{app_exp2}), as shown in Fig.\ \ref{fig_burgers_stability2}. Note the accuracy of the asymptotic approximation along the lower branch even though we only calculated three terms of the asymptotic expansion.

\section{Conclusion} \label{sec6}

This research proposes neural networks for solving, constructing bifurcations, and analyzing the linear stability of nonlinear PDEs. The construction of bifurcations is achieved by incorporating a pseudo-arclength continuation into the network. The method is tested on the one- and two-dimensional Bratu equations and the one-dimensional Burgers equation and is compared to FD.

Compared to FD, NN can have a smaller Jacobian matrix, resulting in reduced memory usage. This advantage is particularly useful for handling higher-dimensional PDEs. Experimental results demonstrate that NN achieves better accuracy and produces more precise bifurcation diagrams compared to FD when using the same grid points. Although NN requires slightly more computational time, this can be mitigated by using fewer grid points. However, in this part of our study, we maintained the same grid settings as FD to ensure a fair comparison. Additionally, the proposed NN can be employed for linear stability analysis by determining the largest eigenvalue of the system.

Note that when we solve the one-dimensional Bratu equation with enforcing the boundary condition, the computational time is around three times faster compared to the time for solving the one dimensional Burgers equation without enforcing the boundary condition. This example indicates that enforcing boundary conditions will be beneficial for reducing the computational time.

For future research, it is important to study the effects of varying the distribution of grid points for NN. Dense grid points can be used in areas of the domain where sharp gradients are likely to occur, while fewer grid points are used in other areas.

Although we employ NN with only two hidden layers and a small number of neurons, the proposed method can be extended to deeper networks with more neurons. This example also demonstrates that different problems \cite{raissi2019physics} require varying numbers of hidden layers and neurons. Specific research focused on determining the optimal number of hidden layers and neurons for solving different PDEs is needed. However, when constructing bifurcation diagrams for PDEs, using a small network is recommended for computational time purposes since the search for solutions needs to be repeated many times 
to obtain complete bifurcation diagrams.

In this research, we utilized pseudo-arclength continuation to trace a fold (i.e., turning point or saddle node) bifurcation. The application of NN to study other types of bifurcations, including global ones, remains to be done. 
Using NN to find several eigenvalues, not only the largest one, is also interesting and relevant. Optimal loss functions in that case remains to be studied.

\section*{Data availability}
No data was used for the research described in the article.

\section*{Declaration of competing interest} 
The authors declare that they have no known competing financial interests or personal relationships that could have appeared to influence the work reported in this paper.

\section*{CRediT authorship contribution statement}

\textbf{Muhammad Luthfi Shahab:} Conceptualization, Software, Formal Analysis, Visualization, Writing - original draft, \textbf{Hadi Susanto:} Conceptualization, Supervision, Writing - review \& editing.

\section*{Acknowledgement}
MLS is supported by a four-year Doctoral Research and Teaching Scholarship (DRTS) from Khalifa University. This work was supported by a Faculty Start-Up Grant (No.\ 8474000351/FSU-2021-011) by Khalifa University. HS also acknowledged a Competitive Internal Research Awards Grant (No.\ 8474000413/CIRA-2021-065) and a Research \& Innovation Grant (No.\ 8474000617/RIG-2023-031) from Khalifa University.

\bibliographystyle{elsarticle-num} 
\bibliography{main}

\appendix
\section{Asymptotic analysis for the eigenvalue problems \eqref{eigenvalueproblem} and \eqref{eigenvalueproblem_burgers}} \label{appendix1}

In this section, we will solve the corresponding eigenvalue problem of the one-dimensional Bratu equation \eqref{eigenvalueproblem} and that of the Burgers equation \eqref{eigenvalueproblem_burgers} with the second boundary condition \eqref{burgers_BC2}. Our results will approximate the largest eigenvalue along the lower stability branch in Figs.\ \ref{fig_bratu_stability2} and \ref{fig_burgers_stability2}. 

The eigenvalue problem \eqref{eigenvalueproblem} with the potential \eqref{analyticalbratu} and \eqref{coshtheta} has a general solution that according to \textsc{Maple} can be written in terms of the Hypergeometric functions. We will instead solve the equation asymptotically. 

Looking at the relation \eqref{coshtheta}, it is natural to take $\omega$ as the small control parameter, in which case the potential of the linear eigenvalue problem can be written as
\begin{align}
	C\exp(u) = 8\omega^2 - 8(1-2x)^2\omega^4 + \cdots. 
\end{align}
It is therefore natural to take the following series expansions
\begin{align}
	\lambda = \lambda_0 +\lambda_1\omega^2 + \lambda_2\omega^4 + \cdots,
	\qquad w = w_0+w_1 \omega^2+w_2 \omega^4+\cdots. \label{app_exp}
\end{align}
Substituting them into \eqref{eigenvalueproblem} and collecting terms of the same power $\mathcal{O}\left(\omega^{2j}\right)$, $j=0,1,2,\dots,$ yield the following equation 
\begin{align}
	\mathcal{L} w_j = f_j, 
\end{align}
where $\mathcal{L} = \lambda_0 - d_x^2$ and $w_j$ must satisfy the Dirichlet boundary conditions $w_j(0)=w_j(1)=0$. Some of the $f_j$s are $f_0=0$, $f_1 = (8-\lambda_1)w_0$, $f_2 = (8-\lambda_1)w_1 - (8(1-2x)^2+\lambda_2)w_0$, etc. Solving the homogeneous equation at $\mathcal{O}(\omega^0)$ is straightforward, i.e., 
\begin{align}
	\lambda_0=-\pi^2,
	\qquad w_0 = \sin(\pi x).
	\label{l0}
\end{align}

The next order equations are inhomogeneous. To be solvable, the inhomogeneous term, i.e., $f_j$, must be orthogonal to the null space of the adjoint of operator $\mathcal{L}$. However, $\mathcal{L}$ is self-adjoint, which means that $w_0$ is a basis of the null space of the adjoint operator. The solvability condition is therefore given by:
\begin{align}
	\int_0^1 f_j\,w_0\,dx=0.
\end{align}
We then obtain the following solutions successively:
\begin{align}
	&\lambda_1=8, \qquad w_1 = 0,\label{l1} \\
	&\lambda_2=\frac{16}{\pi^2}-\frac83, \qquad w_2 = \frac8{\pi^2} \bracketround{x^2-x-\frac1{\pi^2}} \sin(\pi x)-\frac{8x}{3\pi}(2x^2-3x+1)\cos(\pi x),\label{l2}
\end{align}
which, using \eqref{app_exp}, provide the leading order approximation of the critical eigenvalue and corresponding eigenfunction of \eqref{eigenvalueproblem}.

Next, we consider the eigenvalue problem of the Burgers equation \eqref{eigenvalueproblem_burgers} with the second boundary condition \eqref{burgers_BC2}. We will follow the same procedure. In our calculations below, we will set $\nu=1$ without loss of generality, which can be obtained from making the scaling: 
\begin{align}
	u\to\nu u, \qquad t\to t/\nu\label{scale},
\end{align} 
in the governing equation \eqref{burgers}. Our convenient small parameter is then $c$.

Proceeding as before, we note that the eigenvalue problem potential can be expanded as 
\begin{align}
	u(x) = (1-x)c-\frac16(1-x)^3c^2+\cdots.
\end{align}
We then take the following series expansions
\begin{align}
	\lambda = \lambda_0 +\lambda_1c + \lambda_2c^2 + \cdots,
	\qquad w = w_0+w_1 c+w_2 c^2+\cdots. \label{app_exp2}
\end{align}
Substituting them into \eqref{eigenvalueproblem_burgers} and collecting terms of the same power $\mathcal{O}\left(c^{j}\right)$, $j=0,1,2,\dots,$ yield the following equation 
\begin{align}
	\mathcal{L} w_j = f_j, 
\end{align}
where the operator $\mathcal{L}$ is the same as that for the Bratu equation (and hence it is also self-adjoint) and $w_j$ now must satisfy the mixed boundary condition $d_xw_j(0)=w_j(1)=0$. Some of the $f_j$s are $f_0=0$, $f_1 = (1-\lambda_1-(1-x)d_x)w_0$, $f_2 = (1-\lambda_1-(1-x)d_x)w_1 - (\frac12(1-x)^2+\lambda_2-\frac16(1-x)^3d_x)w_0$, etc. Solving the homogeneous equation at $\mathcal{O}(c^0)$ is straightforward, i.e., 
\begin{align}
	\lambda_0=-\frac{\pi^2}4, \qquad w_0 = \cos\left(\frac{\pi x}2 \right).\label{l01}
\end{align}

Following the same procedures gives us the solutions below:
\begin{align}
	&\lambda_1=\frac32, \qquad w_1 = -\frac1\pi(1-x)\sin\left(\frac{\pi x}{2}\right)-\frac14\left(x^2-2x-\frac4{\pi^2}\right)\cos\left(\frac{\pi x}{2}\right),\,\label{l11}\\
	&\lambda_2=-\frac13, \qquad w_2 = \sum_{j=1}^{2}K_j\sin\left(\frac{(2j-1)\pi x}{2}\right)+L_j\cos\left(\frac{(2j-1)\pi x}{2}\right),\label{l21}
\end{align}
for some constants $K_j$ and $L_j$. The leading order approximation of the critical eigenvalue along the lower branch in Fig.\ \ref{fig_burgers_stability2} is given by $\nu\lambda$, where $\lambda$ is the expansion in \eqref{app_exp2} with \eqref{l01}-\eqref{l21}. The multiplication with $\nu$ is due to the time scaling \eqref{scale} above.

\end{document}